\documentclass[a4paper]{amsart}  

\usepackage{amssymb} \usepackage{amscd} \usepackage{times}

\def\zbb{\mathbb{Z}}  
  
  \def\phi{\varphi}
 \def\p1{{\mathbb{P}^1_\zbb}}

\newtheorem{Theorem}{\quad Theorem}[section]

\newcommand{\be} {\begin{equation}}
\newcommand{\ee} {\end{equation}}

\begin{document}

\title{ Blow-up analysis for the solutions of an equation with boundary condition.}

\author{Samy Skander Bahoura}

\address{Department of Mathematics, Pierre and Marie Curie University, 75005, Paris, France.}
              
\email{samybahoura@gmail.com} 

\date{}

\maketitle

\begin{abstract}

Under some conditions we give a blow-up analysis for solutions of an equation with Dirichlet boundary condition.

\end{abstract}

\section{Introduction and Main Result}

\smallskip

We set $ \Delta = -\nabla^i \nabla_i$  the Laplace-Beltrami operator on a connected compact Riemannian manifold with boundary $ (M,g)$  of dimension $ n \geq 3 $ with metric $ g $.

We consider the following equation:

\begin{displaymath}  (P)  \left \{ \begin {split} 
      \Delta u + h u  & = V u^{N-1}, \,\,u >0     \,\, &&\text{in} \!\!&& M, \\
                  u  & = 0  \,\,             && \text{on} \!\!&&\partial M.               
\end {split}\right.
\end{displaymath}

Here, we assume the solutions in the sense of distributions and here also in $ C^{2,\alpha}(\bar M), \alpha >0 $ and the operator $ \Delta +h $ is coercive in $ H^1_0(M) $ with $ h $ a smooth function and $ 0< a \leq V \leq b<+\infty $, $ V $ is a smooth function and $ ||V||_{C^{\alpha}} \leq A $ and $ N=\frac{2n}{n-2}$ the critical Sobolev exponent. 

\smallskip

In  [1-28] we have various  estimates and inequalities of type $ \sup +\inf$ and $ \sup\times \inf $. Here we look to blow-up analysis for regular solutions of the previous equation.

Here we assume:

$$ \sup_M u \cdot \int_M u^{N-1} dV_g \leq C \qquad (*) $$

We have:

$$ ||u||_N \leq ||u||_{\infty}^{\frac{1}{N}}||u||_{N-1}^{1-\frac{1}{N}} \leq |M|^{\frac{1}{N}} ||u||_{\infty} $$

Thus, the condition $ (*) $ is equivalent to the condition:

$$||u||_{\infty}^{\frac{1}{N}}||u||_{N-1}^{1-\frac{1}{N}} \leq C. $$

This condition is an $ L^{\infty}-L^{N-1} $-condition which is intermediate between a $ \sup \times \inf $ condition and $ L^N $ condition, see [1] and [7, 8].

Note that in [7, 8] the condition $(*)$ is satisfied. Thus, this condition is a necessary condition. In [7, 8] it is an exemple where this condition $ (*)$ holds.

\smallskip

Our main results are:

\smallskip

We have the following blow-up analysis on all $ \bar M $:

\smallskip

\begin{Theorem} Assume $(*)$ and $ \sup_M u_i \to +\infty $, then there is a subsequence denoted $ (u_j) $ of $ (P) $ and there are  $ k $ points inside $ M $, $ x_1, \ldots, x_k \in M$, and $ m $ points on $ \partial M $, $ y_1,\ldots, y_m \in \partial M $,  which are concentration points of $ u_j $, and $ u_j \to 0 $ in $ C^2_{loc}(\bar M-\{x_1,\ldots x_k, y_1\ldots, y_m \})$.

\end{Theorem}

\smallskip

To prove this theorem, we use condition $ (*) $ for which there  is a concentration of measures inside $ M $ and on $ \partial M $. As mentionned before, the condition $ (*) $ holds in [7, 8] and thus it is a necessary condition. There is an exemple of blow-up phenomenon which imply the condition $ (*) $. Thus, we must suppose this condition $ (*) $ to know if we have blow-up phenomenon for solutions of the problem $(P)$.

\smallskip

And we have the following compactness result on all $ \bar M $:

\begin{Theorem} We assume $ (*)$. Then there is $ \epsilon_n >0 $ such that if $ b\cdot C < \epsilon_n $, then there is $ c >0 $ such that:

$$ \sup_ M u_i \leq c, \,\, \forall i \in {\mathbb N}. $$

\end{Theorem}

\section{Proofs of the Theorems}

\smallskip

Proofs of Theorems 1.1:

\smallskip

We use two measures, on the locally compact space $ M $, $ \mu_1 $, and on the boundary $ \partial M $, $ \mu_2 $, we use De Giorgi-Nash-Moser iterate scheme, see [1] and the elliptic estimates see [19].

we assume:

$$ \sup_M u_i \cdot \int_M u_i^{N-1} dV_g \leq C \qquad (*) $$

which imply:

\smallskip

$$ \int_M V_i u_i^N dV_g \leq b C $$

Which it is sufficient for the blow-up analysis inside $ M $.

\smallskip

\underbar {Step 1:}

\smallskip

We look to the concentration of measure inside $ M $:

$$ \int_{M} (V_i u_i^N \cdot \phi) dV_g \to \mu_1(\phi) \geq 0, \forall \phi \in C_c^0(M), \phi \geq 0. $$

Here we use the De Giorgi-Nash-Moser iterate scheme, see [1] and the paper [7,8] :

\smallskip

also, we have with the fact that $ \sup_M u_i \to + \infty \Rightarrow \int_M u_i^{N-1} dV_g \to 0 $ and up to a subsequence $ u_i \to 0 $ almost everywhere.

a)

$$ x_0\in M, \mu_1(\{x_0\})< \epsilon_n \leq \frac{\omega_n}{2^n} \Rightarrow \exists \delta >0, \int_{B(x_0,\delta)} V_iu_i^N dV_g  < \epsilon_n \leq \frac{\omega_{n}}{2^n}, $$

$$ x_0\in M, \int_{B(x_0,\delta)} V_iu_i^N dV_g  < \epsilon_n \leq \frac{\omega_{n}}{2^n} \Rightarrow ||u_i||_{L^{\infty}(B(x_0,\delta/2))} \leq C' $$

b) there is a finite number of points $ x_j \in M $ such that:

$$\mu_1(\{x_j\}) \geq \epsilon_n,  \,\, \int_{B(x_j,\delta)} V_iu_i^N dV_g  \geq  \epsilon_n, $$

For $ \delta \to 0 $ small enough, we have:

$$ x_j \in M, \,\, \exists \,\, (x_j^i)_i, \,\, x_j^i \to x_j, \,\, u_i(x_j^i) \to +\infty. $$

and because $ u_i \to 0 $ almost everywhere on $ M $, we have on each compact set $ K \subset M-\{x_j, j=1,\ldots, k\} $: 

$$ \sup_K u_i \to 0. $$

Also we have the convergence to a measure: sum of Dirac measures.

\bigskip

\underbar {Step 2:}

We assume the duality theorem for $ v_i=\sup_M u_i \cdot u_i $, which assert that:

$$ ||\nabla v_i||_{L^q(M)} \leq C_q <+\infty, \forall q,  1\leq q<\frac{n}{n-1}, $$

because $ u_i \to 0 $ on each compact set of $ M $ and the Holder inequality applied for, $ r $,  $ \frac{n}{n-1} > r> q \geq 1 $, we obtain near $ \partial M $ (we prove it in charts of points of the boundary):

$$ ||\nabla v_i||_{r} =o(1),\forall r,  1\leq r<\frac{n}{n-1} . $$

For the duality theorem: 

In fact, we solve $ \Delta z+hz=div f, f\in L^{q'} $ with Dirichlet boundary condition ($ L=\Delta +h $ is coercive, we can solve in Hilbert space), we have: $ z\in L^{\infty} $ with $ ||z||_{L^{\infty}} \leq c_q ||f||_{q'} $, see Gilbarg-Trudinger, [19]. We write:

(we smooth the solutions $ f :\to f_j \in C^{\infty}_c(M), z:\to z_j \in C^{\infty} $, $ f_j\to_{L^{q'}} f $ and thus $ z_j$ is a Cauchy sequence in $L^{\infty}$):

$$ \int_M f\cdot \nabla v_i=\int_M \nabla z\cdot \nabla v_i+ h z v_i=\int_M z v_i u_i^{N-2}, $$

thus,
$$ |\int_M f\cdot \nabla v_i | \leq C||z||_{L^{\infty}} \leq C \cdot c_q ||f||_{q'}, $$

Thus,

$$ ||\nabla v_i||_q \leq c'_q <+\infty .$$

( The functions $ z \in L^2(M) $ uniformly, because we multiply the equation of $ z $ by $ z $ and use Cauchy-Schwarz inequality and the Coercivity of the operator).

Also, we use another type of duality theorem:

We solve (in Hilbert space):

$$ \Delta \tilde z+ h\tilde z= g\in L^{q'/2} $$

We can take $ g\in L^{q'} \subset L^{q'/2} $, because $ q'>n >2 $ to take the norm $ L^2 $.

\smallskip

we have (see Gilbarg-Trudinger, [19]): $ \tilde z\in L^{\infty} $ with $ ||\tilde z||_{L^{\infty}} \leq \tilde c_q ||g||_{q'/2} $. We write:

\smallskip

We smooth the solutions:

$$ \int_M g\cdot  v_i=\int_M \nabla \tilde z\cdot \nabla v_i+ h \tilde z v_i=\int_M \tilde z v_i u_i^{N-2}, $$

thus,
$$ |\int_M g\cdot v_i | \leq C||\tilde z||_{L^{\infty}} \leq C \cdot \tilde c_q ||g||_{q'/2}, $$

Thus,

$$ ||v_i||_{s} \leq \tilde c'_s <+\infty, 1 \leq s <\frac{n}{n-2} .$$

And,

$$ ||v_i||_{q} \leq \bar c'_q <+\infty, 1 \leq q <\frac{n}{n-1} .$$

Also, we can use locally, the Poincar\'e inequality to have a uniform bound of $ \int_M v_i^q, 1\leq q<\frac{n}{n-1} $.(Another fact that this quantity is uniformly bounded).

Let's consider $ \eta $ a cuttof function ( in an open set of a chart), then, we apply the  th. of duality to $ \nabla (v_i\eta) $ in $ M $, but with $\eta $, we are in the open set of the chart:

$$ \Delta (v_i\eta)+h(v_i\eta) = V_i v_i u_i^{N-2} \eta + \nabla v_i\cdot \nabla \eta +v_i \Delta \eta , $$

$$ ||\nabla (v_i\eta)||_q \leq c'_q + C\int_M v_i, $$

but in the chart, we have a Lipschitz domain, we can use the Poincar\'e inequality:

$$ ||v_i\eta||_q \leq c''_q ||\nabla (v_i\eta)||_q \leq \tilde c_q+C\int_M v_i, $$

Thus, in charts with cuttof funtions and partition of unity we obtain:

$$ ||v_i||_q \leq \bar c_q+ C\int_M v_i, \forall q, \,\,1\leq q< \frac{n}{n-1}.$$

By contradiction, if $ \int_M v_i \to +\infty $, then we consider $ w_i=\frac{v_i}{\int_M v_i} $ is bounded in $ \dot H_q^1(M) $ and by the Sobolev embedding, thus converge to a constant $ w\equiv k $, in $ \dot H_r^1(M), 1\leq r <q $, but $ tr(w_i)=0\to tr(w)=k $, this is in contradiction with $ \int_M v_i \leq C||v_i||_r=o(1) \int_M v_i, 1=o(1) $. Thus, $ \int_M v_i \leq C <+\infty, \, \forall i $. Thus:

$$ ||v_i||_q \leq \bar c_q, \forall q, \,\,1\leq q< \frac{n}{n-1}.$$

Also, near the boundary, we have (we prove it around points of the boundary in charts):

$$ ||v_i||_r =o(1), \forall r, \,\,1\leq r< \frac{n}{n-1}.$$

Another fact about the boundedness of $ \int_M v_i $:

Let $ \eta $ solution of:

$$ \Delta \eta + h\eta = 1, \,\, {\rm in }\, M, \, \, {\rm and}  \,\, h=0 \,\,{\rm on}  \,\, \partial M. $$

Multiplying the equation of $ v_i $ by $ h $ and integrate by parts, we obtain:

$$ |\int_M h v_i \eta | \leq C, $$

Multplying  the equation of $ h $ by $ v_i $ and integrate by parts, we obtain:

$$ 0 < \int_M v_i \leq C. $$

\underbar {Step 3:}

\smallskip

We consider a concentration of the measure for: 

$$ \int_M V_i v_i u_i^{N-2} \phi dx \to \mu_3(\phi) \geq 0, \phi\geq 0, \phi\in C^0_c(M) .$$

We have, if we consider: $ \Delta  v_i +(h-V_iu_i^{N-2}) v_i = 0 $ and  $ v_i\in L^q, q>1 $ uniformly, and around $ x_0 $ with $ \mu_1(x_0) < \epsilon_n \Rightarrow u_i \in L^{\infty}(B(x_0,r)), r >0 $ uniformly. thus we apply the Harnack inquality to $ v_i $, we have $ v_i \in L^{\infty}(B(x_0, r')), r > r'>0 $ uniformly. Thus, $ \mu_3(x_0) < \epsilon_n $.

\smallskip

Thus,

$$ \mu_1(x_0) <\epsilon_n \Leftrightarrow \mu_3(x_0) < \epsilon_n. $$

and,

$$ \mu_1(x_j) \geq \epsilon_n >0 \Leftrightarrow \mu_3(x_j) \geq \epsilon_n >0.$$

Thus, 

$$ supp \mu_1= supp \mu_3=\{x_1, \ldots x_k\}. $$

Thus, $ v_i \to v $ in $ W^{1,q}_0(M) $ and locally uniformly on $ M $ with $ v $ solution in the sense of distributions of:

$$ \Delta v+h v= \sum_{j=1}^k \mu_{3,j} \delta_{x_j}, \,\, \mu_{3,j}=\mu_3(x_j) \geq \epsilon_n >0. $$

But, the operator $\Delta +h $ is coercive and thus, we have the existence of the Green function with Dirichlet boundary condition, $ G $ of this operator. Comparing $ v $ and $ G_0(\cdot )=\sum_{j=1}^k \mu_{3,j} G(x_j,\cdot)$, we have: $v= G_0 $ almost everywhere.

Thus,

$$ v_i \to G_0, \,\, {\rm in} \,\, W^{1,q}_0(M) \cap C^2(M-\{x_1,\ldots, x_k\}), \,\, {\rm with } \,\, G_0\in C^2(\bar M-\{x_1,\ldots x_k\}). $$

Thus, as in dimension 2, ($ G_0 $ is smooth up to the boundary), we have: 

$$ ||\nabla (v_i-G_0)||_q =o(1), ||v_i-G_0||_q = o(1), 1\leq q<n/(n-1). $$

Thus, we have the same situation as in dimension 2, except that we have singularities inside $ M $.

\smallskip

Also, we have near the boundary (because  for the tests functions near the boundary we have estimates of type, $ \eta_{\epsilon'}$ with $ ||\nabla \eta_{\epsilon'}||_{\infty} \leq C/(\epsilon')^2 $):

$$ |\Omega_{\epsilon'}| \leq (\epsilon')^n  \Rightarrow ||\nabla v_i||_q+||v_i||_q=O((\epsilon')^n), \,\, i\geq i_0(\epsilon'), \,\, n\geq 3. $$

\underbar {Step 4:}

\smallskip

Now, we look to the case of the boundary:

\smallskip

By the maximum principle and the fact that $ \Delta +h$ is coercive and the Dirichlet boundary condition, we have: $ \partial_{\nu} u_i \geq 0 $.

\smallskip

By the Stokes formula and the duality theorem, see the papers of dimension 2, we have concentration phenomenon on the boundary $ \partial M $ of $ v_i=\sup_M u_i \cdot u_i $:

$$ \int_{\partial M} \partial_{\nu} ( \sup_M u_i \cdot u_i ) d\sigma \leq C'', $$

$$ \int_{\partial M} \partial_{\nu} ( \sup_M u_i \cdot u_i )\cdot \phi d\sigma \to \mu_2(\phi) \geq 0, \forall \phi \in C^0(\partial M), \phi \geq 0. $$

We use the equation with $ v_i $ and we consider a test function $ \eta_{\epsilon} $ by a chart near $ y_0 \in \partial M $, and a Dirichlet problem:

$$ \Delta \tilde \eta_{\epsilon} + h \tilde \eta_{\epsilon} =0 \,\, in \,\, M, \,\, \tilde \eta_{\epsilon} = \eta_{\epsilon} \,\, on \,\,\partial M. $$

For this as in dimension 2, we solve:

$$ \Delta \bar \eta_{\epsilon} +h \bar \eta_{\epsilon} =\Delta \eta_{\epsilon}+h\eta_{\epsilon} \,\, in \,\, M, \,\, \bar \eta_{\epsilon} = 0 \,\, on \,\,\partial M, $$

and we take $ \tilde \eta_{\epsilon}=-\bar \eta_{\epsilon}+\eta_{\epsilon}.$

$$ \Delta (v_i\tilde \eta_{\epsilon})= V_i v_i u_i^{N-2}\tilde \eta_{\epsilon} + \nabla v_i\cdot \nabla \tilde \eta_{\epsilon}, $$

and we integrate.

$$ \int_{M} \Delta (v_i\tilde \eta_{\epsilon})=\int_{\partial M} \partial_{\nu}(v_i\tilde \eta_{\epsilon}) =\int_{\partial M} \partial_{\nu}(v_i)\tilde \eta_{\epsilon} \to \mu_2(\tilde \eta_{\epsilon}) \leq \mu_2(J_{\epsilon}) < \epsilon_n,$$

and near the boundary, we have:

$$ ||\nabla v_i||_{r} =o(1), \,\,||v_i||_r =o(1), \forall r, \,\,1\leq r< \frac{n}{n-1}.$$

Thus,

$$ \sup_M u_i \int_{B(y_0,\delta)} V_iu_i^{N-1} dV_g \leq \int_M V_i v_i u_i^{N-2}\tilde \eta_{\epsilon} < \epsilon_n, $$

\smallskip

In fact we have by the duality theorem in a small neighborhood of $ y_0 \in \partial M $ (and the holder inequality) and the convergence inside $ M $:

$$ y_0\in \partial M, \mu_2(\{y_0\}) <\epsilon_n \Rightarrow \sup_M u_i \int_{B(y_0,\delta)} V_iu_i^{N-1} dV_g  < \epsilon_n \Rightarrow\int_{B(y_0,\delta)} V_iu_i^N dV_g  < \epsilon_n $$

First we consider the equation satisfied by $ v_i $  to use the measure $ \mu_2 $, after we consider the equation satisfied by $ u_i $.

\smallskip

and by the De Giorgi-Nash-Moser iterate scheme applied to $ u_i$ with respect to initial equation, we have:

$$ y_0 \in \partial M, \int_{B(y_0,\delta)} V_iu_i^N dV_g  <\epsilon_n \leq  \frac{\omega_{n}}{2^n} \Rightarrow ||u_i||_{L^{\infty}(B(y_0,\delta/2))} \leq C' $$

and,

There is a finite number of points $ y_j \in \partial M $, such that:

$$ \mu_2(\{y_j\}) \geq \epsilon_n>0. $$

$$\mu_2(\{y_j\}) \geq \epsilon_n >0 \Rightarrow \,\, \forall \delta >0, \,\, \sup_M u_i \cdot \int_{B(y_j,\delta /2)} V_iu_i^{N-1} dV_g=\int_{B(y_j,\delta)} V_i v_i u_i^{N-2} \geq  \frac{\epsilon_n}{2}, $$

For $ \delta \to 0 $ small enough and because $ ||v_i||_1=o(1) $ near the boundary, we have:

$$ y_j\in \partial M, \,\,\exists \,\, (y_j^i)_i, \,\, y_j^i \to y_j, \,\, u_i(y_j^i) \to +\infty. $$

and, for all compact $ \tilde K \subset \bar M-\{x_j, j=1,\ldots k, y_j, j=1,\ldots m\}$

$$ \sup_{\tilde K} u_i \to 0.$$

Also, we have the convergence to a measure on the boundary: sum of Dirac measures on the boundary.

\smallskip

By the Harnack inequality on the boundary see, Gilbarg-Trudinger [19], we have $ v_i \to G_0 $ in $ C^2 (\bar M -\{x_1,\ldots, x_k, y_1,\ldots, y_m \})$, up to the boundary. 

We can see that locally, we have the following asymptotic estimate outside the blow-up points:

$$ u_j(x) \approx \frac{1}{\sup_M u_j}. $$ 

And near the boundary and outside the blow-up points:

$$ u_j(x)\approx \frac{o(1)}{\sup_M u_j}. $$

And,

$$ \sup_M u_j \cdot u_j(\cdot) \to G_0=\sum_{j=1}^k\mu_{3,j} G(x_j,\cdot) $$

in $ \, C^2(\bar M-\{x_1, \ldots, x_k, y_1, \ldots, y_m \}) $ with $ \mu_{3,j}=\mu_3(x_j)\geq \epsilon_n >0.$

\bigskip

\smallskip

Proofs of Theorems 1.2:

\smallskip

In the proof of the theorem 1.1. we have determined $ \epsilon_n >0 $ the mass's reference. By contradiction, if $ \sup_M u_i \to +\infty $, there are concentration of measures on $ M $ and on $ \partial M $ with the total mass $ b\cdot C \geq  \epsilon_n >0 $, which is not possible. Thus, the compactness result.

\end{document}